\documentclass[12pt]{article}
\usepackage{amsfonts}
\usepackage{amsmath}
\usepackage{amssymb}
\def\be {\begin{equation}}
\def\ee {\end{equation}}
\def\ba {\begin{eqnarray}}
\def\ea {\end{eqnarray}}

\begin{document}
\title{ Distribution of Normalized Zero-Sets
of Random Entire Functions }
\author{{Yao Weihong }\\
{\small  Department of Mathematics, Shanghai Jiaotong University,}\\
{\small Shanghai 200240, P. R. China,}\\
{\small E-Mail: \ whyao@sjtu.edu.cn  }}
\date{}
\maketitle \baselineskip 18pt
\begin{center}
\begin{minipage}{120mm}
\textbf{A{\sc \textbf{bstract:}} } This paper is concerned with the
distribution of normalized zero-sets of random entire functions.
 The normalization of the zero-set is performed in the same way as that
of the counting function for an entire function in Nevanlinna
theory. The result generalizes the Shiffman and Zelditch theory on
the distribution of the zeroes of random holomorphic sections of
powers for positive Hermitian holomorphic line bundles from
polynomial functions to entire functions. Our result can also be
viewed as the analogy of Nevanlinna's First Main Theorem in the
theory of the distribution of zero-sets of random entire functions.

\par
\
\par
\textbf{Key words: }Hermitian holomorphic line bundles, random entire functions,  counting function, star operator\\
\end{minipage}
\end{center}

\section{ Introduction}

\bigbreak \, \, \,  In 1955, Hammersley  [SV]  proved that the zeros
of random complex Kac  Polynomials
$$f(z) = \sum_{j=0}^Na_jz^j , \, \, \, \, z \in {\mathbb C}$$
 tend to concentrate on the unit circle
 $S^1 = \{z \mid |z| = 1 \}$
 as
$N\rightarrow \infty$ when the coefficients $a_j$ are independent
complex Gaussian random variables of mean 0 and variance 1:
$$
E(a_j) = 0, \, \,  \, \, \, \,
 E(a_j\overline a_k) = \delta_{jk},  \, \, \, \, \, \, E(a_ja_k) = 0.$$
In recent years, there have been renewed interests in the study of
 distribution of zero-sets of random polynomials.
For example,  Shiffman and Zelditch extended the result to
 random (homogeneous) polynomials of several variables. More generally,
they studied the distribution of zeroes of random holomorphic
sections of positive holomorphic line bundles. These recent
developments were partially motived by the potential application of
the theory to  string/M-theory [DSZ-2004, DSZ-2006a, DSZ-2006b], for
it is observed that the
 supersymmetric vacua (``universes'') of
the string/M-theory may be identified with certain critical points
of a holomorphic section (the ``superpotential'') of a Hermitian
holomorphic line bundle over a complex manifold.

\medbreak  In this paper,
 instead of polynomials, we study the
the distribution of normalized zero-sets of random entire functions.
More precisely,  let $f_1(z),\cdots,f_\ell(z)$  be a finite number
of fixed entire functions.  Let
$$
G_n(z)=\sum_{\nu=0}^n\sum_{j_1=1}^\ell\cdots\sum_{j_\nu=1}^\ell
a_{j_1,\cdots,j_\nu}f_{j_1}(z)\cdots f_{j_\nu}(z)
$$
be a random polynomial, where each coefficient
$a_{j_1,\cdots,j_\nu}$ is an indeterminate which satisfies the
Gaussian distribution
$$
\frac{1}{\pi}\,e^{-|z|^2}
$$
on ${\mathbb C}$.  Note that, by the usual notational convention,
for $\nu=0$ the coefficient $a_{j_1,\cdots,j_\nu}$ is a single
indeterminate with the same Gaussian distribution though we have no
values for $j_1,\cdots,j_\nu$.   In the language of probability
theory,
$$\left(a_{j_1,\cdots,j_\nu}\right)_{0\leq \nu \leq n, 1\leq j_1\leq\ell, \cdots, 1\leq j_{\nu}\leq
\ell}$$ are independent complex Gaussian random variables of mean
$0$ and variance $1$.  Let $N_{\ell,n}$ be the number of elements in
$$\left(a_{j_1,\cdots,j_\nu}\right)_{0\leq \nu \leq n, 1\leq j_1\leq\ell,\cdots,1\leq j_{\nu}\leq
\ell},$$ which is
$$
N_{\ell, n}=1+\ell+\ell^2+\cdots+\ell^n.
$$
Let $a_0$ be the single indeterminate $a_{j_1,\cdots,j_\nu}$ when
$\nu=0$. \medbreak Fix $r>0$,  by Cauchy's integral formula (or the
Poinca\'e-Lelong formula)
$$
\frac{1}{n}\sum_{G_n(z)=0,\atop 0<|z|<r}
\left(\log\frac{r}{\left|z\right|}\right)\delta_z
=\frac{\sqrt{-1}}{n\pi}\left(\log\frac{r}{\left|z\right|}\right)
\partial\bar\partial\log\left|G_n(z)\right|\leqno{(*)}
$$
on the punctured disk $0<\left|z\right|<r$, where $\delta_z$ is the
Dirac delta on ${\mathbb C^l}$ at the point $z$ of ${\mathbb C^l}$.
We now consider the normalized counting divisor ${\mathbf
Z}\left(r,G_n\right)$ of $G_n(z)$ on the punctured disk
$0<\left|z\right|<r$ (in the sense of distribution) which is given
by
$${\mathbf Z}\left(r,G_n\right)=
\frac{1}{n}\sum_{G_n(z)=0,\atop
0<|z|<r}\left(\log\frac{r}{\left|z\right|}\right)\delta_z.
$$
By $(*)$, the expectation ${\mathbf E}\left({\mathbf
Z}\left(r,G_n\right)\right)$ of ${\mathbf Z}\left(r,G_n\right)$ is
equal to
$$
\displaylines{\qquad\int_{\left(a_{j_1,\cdots,j_\nu}\right)\in{\mathbb
C}^{N_{\ell,n}}}
\left(\frac{\sqrt{-1}}{n\pi}\left(\log\frac{r}{\left|z\right|}\right)
\partial\bar\partial\log\left|G_n(z)\right|\right)\cdot\hfill\cr\hfill\cdot\prod_{\left(a_{j_1,\cdots,j_\nu}\right)\in{\mathbb
C}^{N_{\ell,n}}}\left(\frac{1}{\pi}\,e^{-\left|a_{j_1\cdots
j_\nu}\right|^2}\frac{\sqrt{-1}}{2} da_{j_1\cdots j_\nu}\wedge
d\overline{a_{j_1\cdots j_\nu}}\right).\qquad\cr}
$$

 The main result of this paper  is as follows:

\bigskip
{\bf Main Theorem.} {\it Let $C$ be the smooth (possibly non-closed)
curve in ${\mathbb C}$ consisting of all the points $z$ such that
$\left|f(z)\right|=\left(\sum_{j=1}^\ell\left|f_j(z)\right|^2\right)^{\frac{1}{2}}=1$
and
$f^\prime(z)=\left(f_1^\prime(z),...,f_l^\prime(z)\right)\not=0$.
Then the limit of  $ {\mathbf E}\left({\mathbf
Z}\left(r,G_n\right)\right)$ is equal to
$\log\frac{r}{\left|z\right|}$ times the sum of
$$\left|f(z)\right|\Xi\left(\left|f(z)\right|\right)\frac{\sqrt{-1}}{2\pi}\partial\bar\partial\log\left|f(z)\right|$$
and the measure on $C$ defined by the $1$-form
$$
\frac{\sqrt{-1}}{2}\sum_{j=1}^\ell\left(f_j(z)
d\overline{f_j(z)}-\overline{f_j(z)} df_j(z)\right),
$$
where
$$
\Xi\left(x\right)=\left\{\begin{matrix}\frac{2}{x}\ \ {\rm when\ \
}x>1\cr 1\ \ {\rm when\ \ }x=1\hfill\cr 0\ \ {\rm when\ \
}x<1.\hfill\cr\end{matrix}\right.
$$}

\bigskip
When $\ell =1$ and $f_1(z)=z$, our theorem recovers the result of
Hammmersley.

\bigskip
Our main theorem can be viewed as the analogy of Nevanlinna's First
Main Theorem in the theory of the distribution of zero-sets of
random entire functions. It is the first step towards understanding
the theory of distribution of zero-sets of random entire
(meromorphic) functions (we call such theory the ``quantum
Nevanlinna theory''.) This paper sets up a frame work about this
theory, and it would be interesting and meaningful to formula and
prove analogy of Nevanlinna's Second Main Theorem.

\bigbreak The rest of paper is organized as follows. In Section 2,
we present several lemmas, which are needed to prove the main
theorem. In Section 3 we give the proof of the main theorem.

\section{Lemmas }

In this section, we present several lemmas which are needed to prove
the main theorem in the next section. These lemmas are about the
convergence of integrals as distributions.

Firstly, we give the real formula on convergence of integrals as
distributions. See the following Lemma 1.

\bigbreak\noindent{\bf Lemma 1. (On the Convergence of Integrals as
Distributions)} {\it The following formula \be
\lim_{n\to\infty}\frac{1}{n}\int_{x=-\infty}^\infty\left(\frac{d^2}{dx^2}\log\sum_{j=0}^n
e^{jx}\right)\varphi(x)dx=\varphi(0) \label{b13} \ee holds for any
compactly supported continuous function $\varphi$ on ${\mathbb R}$.}

\medbreak\noindent{\bf Proof.}  First check that
$$
\frac{d^2}{dx^2}\log\sum_{j=0}^n e^{jx}\geq 0.
$$
Direct computation yields
$$
\frac{d^2}{dx^2}\log\sum_{j=0}^n e^{jx}=\frac{\left(\sum_{j=0}^n
e^{jx}\right)\left(\sum_{j=0}^n j^2e^{jx}\right)-\left(\sum_{j=0}^n
je^{jx}\right)^2}{\left(\sum_{j=0}^n e^{jx}\right)^2}
$$
which is $\geq 0$ because of the following inequality of
Cauchy-Schwarz
$$
\left(\sum_{j=0}^n a_j b_j\right)^2\leq\left(\sum_{j=0}^n
a_j^2\right)\left(\sum_{j=0}^n b_j^2\right)
$$
with $a_j=e^{\frac{j}{2}x}$ and $b_j=je^{\frac{j}{2}x}$.  To finish
the proof it suffices to check that
$$
\lim_{n\to\infty}\frac{1}{n}\int_{x=-\infty}^\infty\left(\frac{d^2}{dx^2}\log\sum_{j=0}^n
e^{jx}\right)dx=1
$$
together with
$$
\lim_{n\to\infty}\frac{1}{n}\int_{x\in{\mathbb
R},\,\left|x\right|>c}\left(\frac{d^2}{dx^2}\log\sum_{j=0}^n
e^{jx}\right)dx=0\quad{\rm for\ \ }c>0.
$$
Direct computation yields
$$
\sum_{j=0}^n e^{jx}=\frac{1-e^{(n+1)x}}{1-e^x}
$$
and
$$
\displaylines{\frac{1}{n}\frac{d^2}{dx^2}\log\sum_{j=0}^n
e^{jx}=\frac{1}{n}\frac{d^2}{dx^2}
\log\frac{1-e^{(n+1)x}}{1-e^x}\cr=
\frac{1}{n}\frac{d^2}{dx^2}\log\left(1-e^{(n+1)x}\right)-
\frac{1}{n}\frac{d^2}{dx^2}\log\left(1-e^x\right)\cr
=\frac{1}{n}\frac{\left(1-e^{(n+1)x}\right)(n+1)^2e^{(n+1)x}-\left(-(n+1)e^{(n+1)x}\right)^2}
{\left(1-e^{(n+1)x}\right)^2}-\frac{1}{n}\frac{\left(1-e^x\right)e^x-\left(-e^x\right)^2}
{\left(1-e^x\right)^2}\cr}
$$
which approaches to $0$ uniformly on $x\leq-c$ as $n\to\infty$ for
any $c>0$.

\bigbreak For $x\geq c$ with a fixed $c>0$, we can rewrite
$$
\displaylines{\frac{1}{n}\frac{d^2}{dx^2}\log\sum_{j=0}^n
e^{jx}=\frac{1}{n}\frac{d^2}{dx^2}
\log\frac{1-e^{(n+1)x}}{1-e^x}\cr=\frac{1}{n}\frac{d^2}{dx^2}
\log\left(-e^{(n+1)x}\frac{1-e^{-(n+1)x}}{1-e^x}\right)\cr
=\frac{1}{n}\frac{d^2}{dx^2}\left(-(n+1)x\right)+
\frac{1}{n}\frac{d^2}{dx^2}\log\frac{1-e^{-(n+1)x}}{1-e^x}\cr
=\frac{1}{n}\frac{d^2}{dx^2}\log\frac{1-e^{-(n+1)x}}{1-e^x}\cr
=\frac{1}{n}\frac{\left(1-e^{-(n+1)x}\right)(n+1)^2e^{-(n+1)x}-\left((n+1)e^{-(n+1)x}\right)^2}
{\left(1-e^{-(n+1)x}\right)^2}-\frac{1}{n}\frac{\left(1-e^x\right)e^x-\left(-e^x\right)^2}
{\left(1-e^x\right)^2}\cr}
$$
which approaches to $0$ uniformly on $x\geq c$ as $n\to\infty$ for
any $c>0$.

\bigbreak Thus for any fixed $c>0$ we have
$$
\displaylines{\lim_{n\to\infty}\frac{1}{n}\int_{x=-\infty}^\infty\left(\frac{d^2}{dx^2}\log\sum_{j=0}^n
e^{jx}\right)dx=
\lim_{n\to\infty}\frac{1}{n}\int_{x=-c}^c\left(\frac{d^2}{dx^2}\log\sum_{j=0}^n
e^{jx}\right)dx\cr
=\lim_{n\to\infty}\frac{1}{n}\left[\frac{d}{dx}\log\sum_{j=0}^n
e^{jx}\right]_{x=-c}^{x=c}=\lim_{n\to\infty}\frac{1}{n}\left[\frac{d}{dx}\log\frac{1-e^{(n+1)x}}{1-e^x}\right]_{x=-c}^{x=c}\cr
=\lim_{n\to\infty}\frac{1}{n}\left[\frac{-(n+1)e^{(n+1)x}}{1-e^{(n+1)x}}-\frac{-e^x}{1-e^x}\right]_{x=-c}^{x=c}\cr
=\lim_{n\to\infty}\frac{1}{n}\left(\frac{-(n+1)e^{(n+1)c}}{1-e^{(n+1)c}}-\frac{-(n+1)e^{-(n+1)c}}{1-e^{-(n+1)c}}
-\frac{-e^c}{1-e^c}+\frac{-e^{-c}}{1-e^{-c}}\right)\cr
=\lim_{n\to\infty}\frac{1}{n}\,\frac{-(n+1)e^{(n+1)c}}{1-e^{(n+1)c}}\cr
=\lim_{n\to\infty}\frac{1}{n}\,\frac{(n+1)}{1-e^{-(n+1)c}}=1.\cr}
$$
Q.E.D

Now we generalize the result in Lemma 1 from the real domain
${\mathbb R}$ to the complex domain ${\mathbb C}$ as the following
Lemma 2:

\bigbreak\noindent{\bf Lemma 2. (Complex Version of Lemma on the
Convergence of Integrals as Distributions) } {\it The following \be
\lim_{n\to\infty}\frac{1}{n}\int_{\mathbb
C}\left(\frac{\sqrt{-1}}{2\pi}\partial\bar\partial\log\sum_{k=0}^n\left|z\right|^{2k}\right)\varphi(z)
=\frac{1}{2\pi}\int_{\theta=0}^{2\pi}\varphi\left(e^{\sqrt{-1}\theta}\right)d\theta
\label{b17} \ee is true for any compactly supported continuous
function $\varphi$ on ${\mathbb C}$.}

\bigbreak\noindent{\bf Proof.} \  Locally use the holomorphic
coordinate $\log r+\sqrt{-1}\theta$ when $z=re^{\sqrt{-1}\,\theta}$.
Since
$$
\sqrt{-1}\partial\bar\partial
F=\frac{1}{2}\left(\left(\frac{\partial^2}{\partial
x^2}+\frac{\partial^2}{\partial y^2}\right)F\right)\left(dx\wedge
dy\right)
$$
for any function $F$, when we replace $z=x+\sqrt{-1}y$ by $\log
r+\sqrt{-1}\theta$ we get
$$
\sqrt{-1}\partial\bar\partial
F=\frac{1}{2}\left(\left(\frac{\partial^2}{\left(\partial\log
r\right)^2}+\frac{\partial^2}{\partial
\theta^2}\right)F\right)\left(d\log r\wedge d\theta\right).
$$
In the case of
$$F=\log\sum_{k=0}^n\left|z\right|^{2k}=\log\sum_{k=0}^nr^{2k}$$ we
get (when $\xi=2\log r$)
$$
\displaylines{\sqrt{-1}\partial\bar\partial\log\sum_{k=0}^n\left|z\right|^{2k}=\frac{1}{2}\left(\left(\frac{\partial^2}{\left(\partial\log
r\right)^2}\right)\log\sum_{k=0}^nr^{2k}\right)\left(d\log r\wedge
d\theta\right)\cr
=\left(\frac{d^2}{d\xi^2}\log\sum_{k=0}^ne^{k\xi}\right)\left(d\xi\wedge
d\theta\right).\cr}
$$
Thus
$$
\displaylines{\lim_{n\to\infty}\frac{1}{n}\int_{\mathbb
C}\left(\frac{\sqrt{-1}}{2\pi}\partial\bar\partial\log\sum_{k=0}^n\left|z\right|^{2k}\right)\varphi(z)\cr
=\lim_{n\to\infty}\frac{1}{n}\int_{\mathbb
C}\left(\frac{d^2}{d\xi^2}\log\sum_{k=0}^ne^{k\xi}\right)\varphi(e^{\frac{\xi}{2}+\sqrt{-1}\theta})\left(d\xi\wedge
\frac{d\theta}{2\pi}\right)\cr
=\frac{1}{2\pi}\int_{\theta=0}^{2\pi}\varphi\left(e^{\sqrt{-1}\theta}\right)d\theta}
$$
after we integrate with respect to $\xi$ first from $\xi=-\infty$ to
$\infty$ and apply the Lemma on the Convergence of Integrals as
Distributions with $x$ replaced by $\xi$.  Q.E.D.

\bigbreak\noindent {\bf Remark 1. \ } When formulated completely in
terms of distributions the Complex Version of Lemma on the
Convergence of Integrals as Distributions reads
$$
\lim_{n\to\infty}\frac{1}{n}\left(\frac{\sqrt{-1}}{2\pi}\partial\bar\partial\log\sum_{k=0}^n\left|z\right|^{2k}\right)
=\delta_{{}_{S^1}},
$$
where $\delta_{{}_{S^1}}$ means the measure supported on the unit
circle $$S^1=\left\{\,z\in{\mathbb C}\,\Big|\,|z|=1\,\right\}$$
which is the standard measure on $S^1$ of total measure $1$.

\bigbreak\noindent {\bf Remark 2.}  Another way of writing down the
formula in the preceding Lemma is
$$
\lim_{n\to\infty}\frac{1}{n}\left(\frac{\sqrt{-1}}{2\pi}\partial\bar\partial\log\sum_{k=0}^n\left|z\right|^{2k}\right)
=\left.d\left(\frac{1}{2\pi}\,\arg z\right)\right|_{\partial\Delta},
$$
where $\arg z$ is the argument of $z$ and $\Delta$ is the open unit
disk in ${\mathbb C}$ and $\partial\Delta$ is its boundary
$\left\{\left|z\right|=1\right\}$. By using
$$
\displaylines{\partial\bar\partial\Phi(r)=\Phi^\prime(r)\partial\bar\partial
r+\Phi^{\prime\prime}(r)\partial r\wedge\bar\partial r\cr=
\Phi^\prime(r)\left(r\partial\bar\partial\log r+\frac{\partial
r\wedge\bar\partial r}{r}\right)+\Phi^{\prime\prime}(r)\partial
r\wedge\bar\partial r\cr}
$$

with $$\Phi(r)=\frac{1}{n}\log\sum_{k=0}^n r^{2k}$$ when
$r=\left|z\right|$ and $z\in{\mathbb C}$, we can rewrite the limit
as
$$
\lim_{n\to\infty}\left(\frac{1}{n}\left(\frac{1}{r}\,\frac{d}{dr}\log\sum_{k=0}^n
r^{2k}\right)+ \frac{1}{n}\left(\frac{d^2}{dr^2}\log\sum_{k=0}^n
r^{2k}\right)\right)\frac{\sqrt{-1}}{2\pi}\partial r\wedge
\bar\partial r=\left.d\left(\frac{1}{2\pi}\,\arg
z\right)\right|_{\partial\Delta}.
$$

\
\par
\

\medbreak\noindent {\bf Lemma 3.}  \ {\it For
$z=\left(z_1,\cdots,z_\ell\right)\in{\mathbb C}^\ell$ and
$r=\left|z\right|$ and $B_{r_0}=\left\{r<r_0\right\}\subset{\mathbb
C}^\ell$,
$$
\lim_{n\to\infty}\frac{1}{n}\frac{d}{dr}\left(\log\sum_{k=0}^n\left|z\right|^{2k}\right)
\frac{\sqrt{-1}}{2\pi}\partial\bar\partial\log
r=\Xi\left(r_0\right)\frac{\sqrt{-1}}{2\pi}\partial\bar\partial \log
r$$ where
$$
\Xi\left(x\right)=\left\{\begin{matrix}\frac{2}{x}\ \ {\rm when\ \
}x>1\cr 1\ \ {\rm when\ \ }x=1\hfill\cr 0\ \ {\rm when\ \
}x<1.\hfill\cr\end{matrix}\right.
$$
}
 \medbreak\noindent{\bf Proof.}  We distinguish between the
following two cases in our verification.
\begin{itemize}\item[(i)]  When $r_0=1$,
$$\displaylines{\left.\frac{1}{n}\frac{d}{dr}\left(\log\sum_{k=0}^n\left|z\right|^{2k}\right)\right|_{r=r_0}=
\frac{1}{n}\frac{d}{dr}\left(\log\sum_{k=0}^n
r^{2k}\right)\bigg|_{r=r_0}\cr= \left.\frac{1}{n}\frac{\sum_{k=1}^n
2k r^{2k-1}}{\sum_{k=0}^n r^{2k}}\right|_{r=r_0}
=\frac{1}{n}\frac{\sum_{k=1}^n 2k}{\sum_{k=0}^n 1}= \frac{1}{n}
\frac{n(n+1)}{n+1}=1.\cr}
$$
\item[(ii)] When $r_0\not=1$,
$$\displaylines{\left.\frac{1}{n}\frac{d}{dr}\left(\log\sum_{k=0}^n\left|z\right|^{2k}\right)\right|_{r=r_0}=
\left.\frac{1}{n}\frac{d}{dr}\left(\log\sum_{k=0}^n
r^{2k}\right)\right|_{r=r_0}\cr=
\left.\frac{1}{n}\frac{d}{dr}\left(\log\frac{1-r^{2(n+1)}}{1-r^2}\right)\right|_{r=r_0}
=
\left.\frac{1}{n}\frac{d}{dr}\left(\log\left(1-r^{2(n+1)}\right)-\log\left(1-r^2\right)\right)\right|_{r=r_0}\cr
=
\left.\frac{1}{n}\left(\frac{-(2n+2)r^{2n+1}}{1-r^{2(n+1)}}+\frac{2r}{1-r^2}\right)\right|_{r=r_0},\cr}
$$
whose limit is $ \frac{2}{r_0}$ as $n\to\infty$ when $r_0>1$ and
whose limit is $0$ when $0\leq r_0<1$.
\end{itemize}

\bigbreak\noindent {\bf Lemma 4. \ (Differential of Argument in a
Prescribed Complex Line)}.

{\it  The differential of the argument on the complex line in
${\mathbb C}^\ell$ through the point
$\left(z_1^{(0)},\cdots,z_\ell^{(0)}\right)$ and the origin is given
by
$$
\frac{\sqrt{-1}}{2}\left(\frac{d\overline{\left<z,z^{(0)}\right>}}{\overline{\left<z,z^{(0)}\right>}}-\frac{d\left<z,z^{(0)}\right>}{\left<z,z^{(0)}\right>}\right).
$$
In particular, the value at $z\in{\mathbb C}^\ell$ of the
differential of the argument on the complex line in ${\mathbb
C}^\ell$ through the point $z\in{\mathbb C}^\ell$ and the origin is
equal to
$$
\frac{1}{\left|z\right|^2}\frac{\sqrt{-1}}{2}\sum_{j=1}^\ell\left(z_j
d\bar z_j- \bar z_jdz_j\right).
$$
}

\medbreak\noindent{\bf Proof.}   Consider the circle consisting of
all $z\in{\mathbb C}^\ell$ such that
$\left|z\right|=\left|z^{(0)}\right|$ and $z=\lambda z^{(0)}$ for
some $\lambda\in{\mathbb C}$, which means
$\left<z,z^{(0)}\right>=\lambda\left|z^{(0)}\right|^2$, or
$z=\frac{\left<z,z^{(0)}\right>}{\left|z^{(0)}\right|^2} z_0$.  The
differential of the argument ${\rm arg}\lambda$ is given by
$\frac{1}{2\sqrt{-1}}$ times the differential of
$\log\frac{\lambda}{\bar\lambda}$, which means
$$
\frac{\sqrt{-1}}{2}\left(\frac{d\overline{\left<z,z^{(0)}\right>}}{\overline{\left<z,z^{(0)}\right>}}-\frac{d\left<z,z^{(0)}\right>}{\left<z,z^{(0)}\right>}\right).
$$
When we do the evaluation at the point $z^{(0)}$, we have
$$
\displaylines{ d\overline{\left<z,z^{(0)}\right>}=\sum_{j=1}^\ell
z^{(0)}_j d\bar z_j,\cr
d\left<z,z^{(0)}\right>=\sum_{j=1}^\ell\overline{z^{(0)}_j}dz_j,\cr
\left<z,z^{(0)}\right>=\left|z^{(0)}\right|^2.\cr}
$$
Q.E.D.

\bigbreak\noindent {\bf Lemma 5 ( Lemma with Pullback by Holomorphic
Function).} \ {\it Let $\Omega$ be a connected open subset of
${\mathbb C}$ and $f:\Omega\to{\mathbb C}$ be a nonconstant
holomorphic function on $\Omega$. Let $C$ be the smooth (possibly
non-closed) curve in $\Omega$ consisting of all the points $z$ of
$\Omega$ such that $\left|f(z)\right|=1$ and $f^\prime(z)\not=0$.
Then
$$
\lim_{n\to\infty}\frac{1}{n}\left(\frac{\sqrt{-1}}{2\pi}\partial\bar\partial\log\sum_{k=0}^n\left|f(z)\right|^{2k}\right)
=\left.d\left(\frac{1}{2\pi}\,\arg f\right)\right|_C
$$
on $\Omega$.}

\medbreak\noindent{\bf Proof.}  The formula follows from pulling
back the formula in Remark 2 by $f$.

\section{The Proof of Main Theorem }

\bigbreak\noindent {\bf Proposition 1 (Complex Version of Lemma on
the Convergence of Integrals as Distributions).}\ {\it
$$
\displaylines{\lim_{n\to\infty}\frac{1}{n}\left(\frac{\sqrt{-1}}{2\pi}\partial\bar\partial\log\sum_{k=0}^n\left|z\right|^{2k}\right)\cr
=r\Xi\left(r\right)\frac{\sqrt{-1}}{2\pi}\partial\bar\partial\log r+
\left[\delta_{{}_{S^{2l-1}}}\right]\wedge
\frac{1}{r^2}\frac{\sqrt{-1}}{2}\sum_{j=1}^\ell\left(z_j d\bar
z_j-\bar z_j dz_j\right),\cr}
$$
where $\left[\delta_{{}_{S^{2l-1}}}\right]$ denotes the $1$-current
on ${\mathbb C}^\ell$ defined by integration over
$$S^{2l-1}=\left\{\,z\in{\mathbb C^l}\,\left|\,|z|=\left(\sum_{j=1}^\ell\left|z_j\right|^2\right)^{\frac{1}{2}}=1\,\right.\right\},$$
and $r=|z|$.}

\bigbreak\noindent{\bf Proof.}  Before giving the proof we would
like to remark that the case of $\ell=1$ is simply Lemma 5.

\medbreak For the proof of the case of a general $\ell$, consider
the natural projection $\pi:{\mathbb
C}^\ell-\left\{0\right\}\to{\mathbb P}_{\ell-1}$ defined by
$\left(z_1,\cdots,z_\ell\right)\to\left[z_1,\cdots,z_\ell\right]$.
Fix $r_0>0$.  Let $P_{r_0}\in{\mathbb C}^\ell-\left\{0\right\}$ be
defined by $z_1=r_0$ and $z_2=\cdots=z_\ell=0$ and let
$Q_0=\pi\left(P_{r_0}\right)$ which is independent of $r_0>0$.
Consider the inhomogeneous coordinate system
$$\left(\frac{z_2}{z_1},\cdots,\frac{z_\ell}{z_1}\right)$$ of
${\mathbb P}_{\ell-1}$ centered at $Q_0$.

\medbreak We choose a special coordinate system of differential
$1$-forms for points of ${\mathbb C}^\ell-\left\{0\right\}$.  At the
point $P_{r_0}$, we choose the differentials $d\left|z_1\right|$,
$d\left(\arg z_1\right)$ of the polar coordinates of $z_1$ and the
differentials
$$d\left(\pi^*\left(\frac{z_2}{z_1}\right)\right),\cdots,d\left(\pi^*\left(\frac{z_\ell}{z_1}\right)\right),
d\overline{\left(\pi^*\left(\frac{z_2}{z_1}\right)\right)},\cdots,d\overline{\left(\pi^*\left(\frac{z_\ell}{z_1}\right)\right)}$$
of the pullback by $\pi$ of the inhomogeneous coordinate system
$$\left(\frac{z_2}{z_1},\cdots,\frac{z_\ell}{z_1}\right)$$ of
${\mathbb P}_{\ell-1}$ centered at $Q_0$ and their complex
conjugates.

\medbreak At $P_{r_0}$, we have
$$
\frac{1}{r^2}\frac{\sqrt{-1}}{2}\sum_{j=1}^\ell\left(z_j d\bar
z_j-\bar z_j dz_j\right)=\frac{1}{r_0}\,d\left(\arg
z_1\right)\leqno{(\sharp)}
$$
either by direct computation or from Lemma 4.

\medbreak Let $U(\ell)$ act on ${\mathbb C}^\ell$ and let
$U\left(\ell-1\right)$ be the stability subgroup of
$U\left(\ell\right)$ for the complex line $z_2=\cdots=z_\ell=0$ in
${\mathbb C}^\ell$.  We take a local smooth submanifold $M$ of real
dimension $2\ell-1$ in $U\left(\ell\right)$ which contains the
identity element of $U\left(\ell\right)$ and is transversal to
$U\left(\ell-1\right)$.   By replacing $r_0$ by a varying $r_1$
which is very close to $r_0$ and by using the action of a varying
element $g$ of $M$, we get a smooth basis of $1$-forms
$\omega_1,\cdots,\omega_{2\ell}$ on an open neighborhood $W$ of
$P_{r_0}$ from the basis of $1$-forms
$$\displaylines{\left(d\left|z_1\right|\right)_{P_{r_0}},  \left(d\left(\arg z_1\right)\right)_{P_{r_0}}, \left(d\left(\pi^*\left(\frac{z_2}{z_1}\right)\right)\right)_{P_{r_0}},\cdots,\left(d\left(\pi^*\left(\frac{z_\ell}{z_1}\right)\right)\right)_{P_{r_0}},\cr \left(d\overline{\left(\pi^*\left(\frac{z_2}{z_1}\right)\right)}\right)_{P_{r_0}},\cdots,\left(d\overline{\left(\pi^*\left(\frac{z_\ell}{z_1}\right)\right)}\right)_{P_{r_0}}\cr}$$
at $P_0$.

\medbreak Now use
$$
\displaylines{\partial\bar\partial\Phi(r)=\Phi^\prime(r)\partial\bar\partial
r+\Phi^{\prime\prime}(r)\partial r\wedge\bar\partial r\cr=
\Phi^\prime(r)\left(r\partial\bar\partial\log r+\frac{\partial
r\wedge\bar\partial r}{r}\right)+\Phi^{\prime\prime}(r)\partial
r\wedge\bar\partial r\cr}
$$
with $$\Phi(r)=\frac{1}{n}\log\sum_{k=0}^n r^{2k}$$ to get
$$
\displaylines{ \frac{1}{n}\partial\bar\partial\log\sum_{k=0}^n
r^{2k}=\frac{1}{n}\left(\frac{d}{dr}\log\sum_{k=0}^n
r^{2k}\right)r\partial\bar\partial\log
r\hfill\cr\hfill+\left(\frac{1}{n}\left(\frac{d}{dr}\log\sum_{k=0}^n
r^{2k}\right)\frac{\partial r\wedge\bar\partial
r}{r}+\frac{1}{n}\left(\frac{d^2}{dr^2}\log\sum_{k=0}^n
r^{2k}\right)\partial r\wedge \bar\partial r\right).\cr}
$$
Since by Lemma 3
$$
\lim_{n\to\infty}\frac{1}{n}\left(\frac{d}{dr}\log\sum_{k=0}^n
r^{2k}\right)r\partial\bar\partial\log
r=r\Xi(r)\partial\bar\partial\log r,
$$
to finish the proof it suffices to verify that
$$
\displaylines{(\natural)\qquad\qquad
\lim_{n\to\infty}\left(\frac{1}{n}\left(\frac{1}{r}\,\frac{d}{dr}\log\sum_{k=0}^n
r^{2k}\right)+ \frac{1}{n}\left(\frac{d^2}{dr^2}\log\sum_{k=0}^n
r^{2k}\right)\right)\frac{\sqrt{-1}}{2\pi}\partial r\wedge
\bar\partial r\cr\hfill=\left[\delta_{{}_{S^{2l-1}}}\right]\wedge
\frac{1}{r^2}\frac{\sqrt{-1}}{2}\sum_{j=1}^\ell\left(z_j d\bar
z_j-\bar z_j dz_j\right).}
$$
Since $S^{2\ell-1}$ is defined by $r=1$, it follows that the
$1$-form $\left[\delta_{{}_{S^{2l-1}}}\right]$ defined by
integration over $S^{2\ell-1}$ is equal to $dr$ times a generalized
function.  It follows from $(\sharp)$ that the right-hand side
$$\left[\delta_{{}_{S^{2l-1}}}\right]\wedge
\frac{1}{r^2}\frac{\sqrt{-1}}{2}\sum_{j=1}^\ell\left(z_j d\bar
z_j-\bar z_j dz_j\right)$$ of $(\natural)$ is equal to a generalized
function times $d\left|z_1\right|\wedge d\left({\rm
arg\,}z_1\right)$.  Since at $P_{r_0}$ we have
$$
\frac{\sqrt{-1}}{2\pi}\partial r\wedge\bar\partial
r=\frac{r_0}{\pi}d\left|z_1\right|\wedge d\left({\rm
arg\,}z_1\right)
$$
and since both sides of $(\natural)$ are invariant under $U(\ell)$,
in order to verify $(\natural)$ it suffices to use a test
$2\ell$-form $\varphi$ supported on $W$ whose coefficients of the
terms containing either $\omega_1$ or $\omega_2$ as a factor are
zero when expressed in terms of $\omega_1,\cdots,\omega_{2\ell}$. In
other words, we need only use a test $2\ell$-form $\varphi$ which is
equal to a function supported on $W$ times the pullback by $\pi$ of
the Fubini-Study volume form of ${\mathbb P}_{\ell-1}$.  When we do
the integration of both sides of $(\natural)$ against $\varphi$, we
can use Fubini's theorem to integrate along the fibers of $\pi$
first.  The integration along each fiber on both sides of
$(\natural)$ gives the same result because of Remark 2.  Thus
$(\natural)$ is verified.  This finishes our proof.  Q.E.D.

\bigbreak

 Now we are ready to prove the main theorem. By
$(*)$, the expectation ${\mathbf E}\left({\mathbf
Z}\left(r,G_n\right)\right)$ of ${\mathbf Z}\left(r,G_n\right)$ is
equal to
$$
\displaylines{\qquad\int_{\left(a_{j_1,\cdots,j_\nu}\right)\in{\mathbb
C}^{N_{\ell,n}}}
\left(\frac{\sqrt{-1}}{n\pi}\left(\log\frac{r}{\left|z\right|}\right)
\partial\bar\partial\log\left|G_n(z)\right|\right)\cdot\hfill\cr\hfill\cdot\prod_{\left(a_{j_1,\cdots,j_\nu}\right)\in{\mathbb
C}^{N_{\ell,n}}}\left(\frac{1}{\pi}\,e^{-\left|a_{j_1\cdots
j_\nu}\right|^2}\frac{\sqrt{-1}}{2} da_{j_1\cdots j_\nu}\wedge
d\overline{a_{j_1\cdots j_\nu}}\right).\qquad\cr}
$$
We introduce two column vectors
$$
\vec{\mathbf a}=\left[a_{j_1,\cdots,j_\nu}\right]_{0\leq\nu\leq n,
1\leq j_1\leq\ell,\cdots,1\leq j_\nu\leq \ell}
$$ and
$$\vec {\mathbf
v}(z)=\left[f_{j_1}(z)\cdots f_{j_\nu}(z)\right]_{0\leq\nu\leq n,
1\leq j_1\leq\ell,\cdots,1\leq j_\nu\leq \ell}
$$
of $N_{\ell,n}$ components each. Here we set $f_0(z)=1$. Then
$G_n(z)$ is equal to the inner product
$$
\left<\vec{\mathbf a},\,\vec{\mathbf v}(z)\right>=\sum_{0\leq\nu\leq
n, 1\leq j_1\leq\ell,\cdots,1\leq j_\nu\leq
\ell}a_{j_1,\cdots,j_\nu}f_{j_1}(z)\cdots f_{j_\nu}(z)
$$
of the two $N_{\ell,n}$-vectors $\vec{\mathbf a}$ and $\vec{\mathbf
v}(z)$. The length of the $N_{\ell,n}$-vector $\vec{\mathbf v}(z)$
is given by
$$\left\|\vec{\mathbf
v}(z)\right\|=\left(\sum_{0\leq\nu\leq n, 1\leq
j_1\leq\ell,\cdots,1\leq j_\nu\leq
\ell}|f_{j_1}(z)|^2\cdots|f_{j_\nu}(z)|^2\right)^{\frac{1}{2}}.$$
Introduce the unit $N_{\ell,n}$-vector
$$
\displaylines{\vec{\mathbf u}(z)=\frac{1}{\left\|\vec{\mathbf
v}(z)\right\|}\,{\mathbf v}(z)\cr=\frac{1}{\left(\sum_{0\leq\nu\leq
n, 1\leq j_1\leq\ell,\cdots,1\leq j_\nu\leq
\ell}|f_{j_1}(z)|^2\cdots|f_{j_\nu}(z)|^2\right)^{\frac{1}{2}}}\left[f_{j_1}(z)\cdots
f_{j_\nu}(z)\right]_{0\leq\nu\leq n, 1\leq j_1\leq\ell,\cdots,1\leq
j_\nu\leq \ell}\cr}
$$
in the same direction as $\vec{\mathbf v}(z)$.  Then
$$
\displaylines{\log\left|G_n(z)\right|=\log\left|\left<\vec{\mathbf
a},\,\vec{\mathbf v}(z)\right>\right|=\log\left|\left<\vec{\mathbf
a},\,\left\|\vec{\mathbf v}(z)\right\|\vec{\mathbf
u}(z)\right>\right|\cr =\log\left\|\vec{\mathbf
v}(z)\right\|+\log\left|\left<\vec{\mathbf a},\,\vec{\mathbf
u}(z)\right>\right|.\cr}
$$
Now ${\mathbf E}\left({\mathbf Z}\left(r,G_n\right)\right)$ is equal
to
$$
\displaylines{\qquad\qquad\int_{\left(a_{j_1,\cdots,j_\nu}\right)\in{\mathbb
C}^{N_{\ell,n}}}
\left(\frac{\sqrt{-1}}{n\pi}\left(\log\frac{r}{\left|z\right|}\right)\partial\bar\partial\left(\log\left\|\vec{\mathbf
v}(z)\right\|+\log\left|\left<\vec{\mathbf a},\,\vec{\mathbf
u}(z)\right>\right|\right)\right)\cdot\hfill\cr\hfill\cdot\prod_{\left(a_{j_1,\cdots,j_\nu}\right)\in{\mathbb
C}^{N_{\ell,n}}}\left(\frac{1}{\pi}\,e^{-\left|a_{j_1\cdots
j_\nu}\right|^2}\frac{\sqrt{-1}}{2} da_{j_1\cdots j_\nu}\wedge
d\overline{a_{j_1\cdots j_\nu}}\right).\qquad\cr}
$$
Let $ \vec{{\mathbf e}}_0$ be the $N_{\ell,n}$-vector
$$\left(e_{j_1,\cdots,j_\nu}\right)_{0\leq\nu\leq n, 1\leq j_1\leq\ell,\cdots,1\leq j_\nu\leq
\ell}$$ whose only nonzero component is $e_0=1$.  Here comes the key
point of the whole argument. For fixed $z$, we integrate
$$
\displaylines{\qquad\qquad\int_{\left(a_{j_1,\cdots,j_\nu}\right)\in{\mathbb
C}^{N_{\ell,n}}}\left(\frac{\sqrt{-1}}{n\pi}\partial\bar\partial\log\left|\left<\vec{\mathbf
a},\,\vec{\mathbf
u}(z)\right>\right|\right)\cdot\hfill\cr\hfill\cdot\prod_{\left(a_{j_1,\cdots,j_\nu}\right)\in{\mathbb
C}^{N_{\ell,n}}}\left(\frac{1}{\pi}\,e^{-\left|a_{j_1\cdots
j_n}\right|^2}\frac{\sqrt{-1}}{2} da_{j_1\cdots j_n}\wedge
d\overline{a_{j_1\cdots j_n}}\right).\qquad\cr
=\int_{\left(a_{j_1,\cdots,j_\nu}\right)\in{\mathbb
C}^{N_{\ell,n}}}\left(\frac{\sqrt{-1}}{n\pi}\partial\bar\partial\log\left|\left<\vec{\mathbf
a},\,\vec{\mathbf
u}(z)\right>\right|\right)\frac{1}{\pi^{N_{\ell,n}}}\,
e^{-\left\|\vec{\mathbf a}\right\|^2}\cr
=\int_{\left(a_{j_1,\cdots,j_\nu}\right)\in{\mathbb C}^{N_{\ell,n}}}
\left(\frac{\sqrt{-1}}{n\pi}\partial\bar\partial\log\left|\left<\vec{\mathbf
a},\,\vec{{\mathbf e}}_0\right>\right|\right)
\frac{1}{\pi^{N_{\ell,n}}}\, e^{-\left\|\vec{\mathbf
a}\right\|^2}\cr =\int_{\left(a_{j_1,\cdots,j_\nu}\right)\in{\mathbb
C}^{N_{\ell,n}}}\left(\frac{\sqrt{-1}}{n\pi}\partial\bar\partial\log\left|a_0
\right|\right)\frac{1}{\pi^{N_{\ell,n}}}\, e^{-\left\|\vec{\mathbf
a}\right\|^2}\cr=\frac{\sqrt{-1}}{n\pi}\partial\bar\partial\int_{a_0\in{\mathbb
C}}\left(\log\left|a_0\right|\right)\frac{1}{\pi}e^{-\left|a_0\right|^2}\cr}
$$
which is equal to
$$\frac{\sqrt{-1}}{n\pi}\partial\bar\partial A=0\quad{\rm with\ \ }A=\int_{a_0\in{\mathbb
C}}\left(\log\left|a_0\right|\right)\frac{1}{\pi}e^{-\left|a_0\right|^2},$$
because $A$ is a constant.  Note that the equality
$$
\displaylines{\int_{\left(a_{j_1,\cdots,j_\nu}\right)\in{\mathbb
C}^{N_{\ell,n}}}\left(\frac{\sqrt{-1}}{n\pi}\partial\bar\partial\log\left|\left<\vec{\mathbf
a},\,\vec{\mathbf
u}(z)\right>\right|\right)\frac{1}{\pi^{N_{\ell,n}}}\,
e^{-\left\|\vec{\mathbf a}\right\|^2}\cr
=\int_{\left(a_{j_1,\cdots,j_\nu}\right)\in{\mathbb C}^{N_{\ell,n}}}
\left(\frac{\sqrt{-1}}{n\pi}\partial\bar\partial\log\left|\left<\vec{\mathbf
a},\,\vec{{\mathbf e}}_0\right>\right|\right)
\frac{1}{\pi^{N_{\ell,n}}}\, e^{-\left\|\vec{\mathbf
a}\right\|^2}\cr}
$$
in the above string of equalities comes from the fact that for any
fixed $z\in{\mathbb C}$ some unitary transformation of ${\mathbb
C}^{N_{\ell,n}}$ (which may depend on $z$) maps $\vec{\mathbf u}(z)$
to $\vec{{\mathbf e}}_0$ and that $e^{-\left\|\vec{\mathbf
a}\right\|^2}$ is unchanged under any unitary transformation acting
on $\vec{\mathbf a}$.  Thus the limit of ${\mathbf E}\left({\mathbf
Z}\left(r,G_n\right)\right)$ as $n\to\infty$ is equal to
$$
\lim_{n\to\infty}\int_{\left(a_{j_1,\cdots,j_\nu}\right)\in{\mathbb
C}^{N_{\ell,n}}}\left(\frac{\sqrt{-1}}{n\pi}\left(\log\frac{r}{\left|z\right|}\right)\partial\bar\partial\log\left\|\vec{\mathbf
v}(z)\right\|\right)\frac{1}{\pi^{N_{\ell,n}}}\,
e^{-\left\|\vec{\mathbf a}\right\|^2},$$ which after integration
over $$\left(a_{j_1,\cdots,j_\nu}\right)_{0\leq\nu\leq n, 1\leq
j_1\leq\ell,\cdots,1\leq j_\nu\leq \ell}$$ is simply equal to
$$
\displaylines{\lim_{n\to\infty}\frac{\sqrt{-1}}{n\pi}\left(\log\frac{r}{\left|z\right|}\right)
\partial\bar\partial\log\left\|
\vec{\mathbf v}(z)\right\|\cr
=\lim_{n\to\infty}\frac{1}{n}\left(\log\frac{r}{\left|z\right|}\right)\left(\frac{\sqrt{-1}}{2\pi}\partial\bar\partial\log\sum_{k=0}^n\left(\sum_{j=1}^\ell
\left|f_j(z)\right|^2\right)^k\right)\cr}
$$
which is equal, by Proposition 1 to $\log\frac{r}{\left|z\right|}$
times the pullback by $f$ of
$$\left|w\right|\Xi\left(\left|w\right|\right)\frac{\sqrt{-1}}{2\pi}\partial\bar\partial\log\left|w\right|+
\left[\delta_{{}_{S^{2l-1}}}\right]\wedge
\frac{1}{\left|w\right|^2}\frac{\sqrt{-1}}{2}\sum_{j=1}^\ell\left(w_j
d\bar w_j-\bar w_j dw_j\right),
$$
where $w\in{\mathbb C}^\ell=\left(w_1,\cdots,w_\ell\right)$ is
variable in the target space of the map
$f=\left(f_1,\cdots,f_\ell\right):{\mathbb C}\to{\mathbb C}^\ell$.
This finishes the proof.

\bigbreak  Now let's consider some spacial cases.

\bigbreak\noindent {\bf Case 1.} \

  Let $\ell=1$, and
$f_1=z$, then $G_n(z)$ become a random polynomial (Kac Polynomial)
$$
G_n(z)=\sum_{\nu=0}^n a_{\nu}z^\nu,
$$
where coefficient $a_{\nu}$ for
 $0<\nu\leq n$ is an
indeterminate which satisfies the Gaussian distribution
$$
\frac{1}{\pi}\,e^{-|z|^2}
$$
on ${\mathbb C}$.   Note that  the total number of indeterminacies

\[ \alpha =
\left(\begin{array}{clr}
a_0  \\
a_1  \\
\vdots  \\
a_n
\end{array}\right)
\]

 is equal to \be N=n+1. \label{b10} \ee

Hence our main theorem implies the following Corollary 1.

\medbreak\noindent{\bf Corollary 1.} \  {\it Let
 $C_1$ be  $\left|z\right|=1$. Then the limit of $ {\mathbf
E}\left({\mathbf Z}\left(r,G_n\right)\right)$ is equal to
$\log\frac{r}{\left|z\right|}$ times the measure on $C_1$ defined by
the $1$-form
$$
\frac{\sqrt{-1}}{2}\left(z d\bar{z}-\bar{z} dz\right).
$$
}  This recovers the result of Hammmersley.

\bigbreak\noindent {\bf Case 2.} \

By restricting to each line and then integrate, our theorem can
actually be extended to $f: {\mathbb C}^\ell \rightarrow {\mathbb
C}^\ell$. By taking $f_1(z)=z_1, \dots, f_{\ell}(z)=z_{\ell}$, then
we get the result given by Shiffman and Zelditch [SZ-1999].

\medbreak\noindent{\bf Remark 3.} \ In the context described above,
we can also investigate in addition the effect of the perturbation
by small functions.

\medbreak\noindent{\bf Open Question.} \ Can  our main theorem be
generalized to the case of meromorphic functions?

\section{Acknowledgments}

\bigbreak The author would like to thank Professors Yum-Tong Siu and
Min Ru for discussions and suggestions about this article.

\section{ References}

 \bigbreak

\medbreak\noindent[DSZ-2004] Michael Douglas, Bernard Shiffman, and
Steve Zelditch, Critical points and supersymmetric vacua. I. {\it
Comm. Math. Phys.} 252 (2004), no. 1-3, 325--358.

\medbreak\noindent[DSZ-2006a] Michael R. Douglas, Bernard Shiffman,
and Steve Zelditch, Critical points and supersymmetric vacua. II.
Asymptotic and extremal metrics. {\it J. Differential Geom.} 72
(2006), no. 3, 381--427.

\medbreak\noindent[DSZ-2006b] Michael Douglas, Bernard Shiffman, and
Steve Zelditch, Critical points and supersymmetric vacua. III.
String/M models. {\it Comm. Math. Phys.} 265 (2006), no. 3,
617--671.

\medbreak\noindent[SV] L.A. Shepp and R.J. Vanderbei, The complex
zeros of random polynomials. {\it Trans. Am. Math. Soc.}
\textbf{347} (1995), 4365--4384.

\medbreak\noindent[SZ-1999] Bernard Shiffman and Steve Zelditch,
Distribution of zeros of random and quantum chaotic sections of
positive line bundles. {\it Comm. Math. Phys.} 200 (1999), no. 3,
661--683.

\end{document}